\documentclass[10pt, leqno]{amsart}

\usepackage{fullpage}
\usepackage{amsmath,amssymb,amsthm, epsfig}

\usepackage{hyperref}

\usepackage{amsmath}

\usepackage[pdftex]{color}

\usepackage{color}

\usepackage{ulem}

\usepackage{epsfig}

\usepackage{mathrsfs}

\usepackage{wasysym}

\usepackage{stackrel}

\usepackage{ifthen}

\usepackage{color, xcolor}

\def \N {\mathbb{N}}
\def \R {\mathbb{R}}

\newcommand{\intav}[1]{\mathchoice {\mathop{\vrule width 6pt height 3 pt depth  -2.5pt
\kern -8pt \intop}\nolimits_{\kern -6pt#1}} {\mathop{\vrule width
5pt height 3  pt depth -2.6pt \kern -6pt \intop}\nolimits_{#1}}
{\mathop{\vrule width 5pt height 3 pt depth -2.6pt \kern -6pt
\intop}\nolimits_{#1}} {\mathop{\vrule width 5pt height 3 pt depth
-2.6pt \kern -6pt \intop}\nolimits_{#1}}}

\newtheorem{theorem}{Theorem}[section]

\newtheorem{lemma}[theorem]{Lemma}

\newtheorem{proposition}[theorem]{Proposition}

\theoremstyle{definition}

\newtheorem{definition}[theorem]{Definition}

\theoremstyle{remark}

\newtheorem{remark}[theorem]{Remark}

\numberwithin{equation}{section}

\allowdisplaybreaks

\def\Xint#1{\mathchoice
{\XXint\displaystyle\textstyle{#1}}%
{\XXint\textstyle\scriptstyle{#1}}%
{\XXint\scriptstyle\scriptscriptstyle{#1}}%
{\XXint\scriptscriptstyle\scriptscriptstyle{#1}}%
\!\int}
\def\XXint#1#2#3{{\setbox0=\hbox{$#1{#2#3}{\int}$ }
\vcenter{\hbox{$#2#3$ }}\kern-.6\wd0}}

\def\dashint{\Xint-}

\begin{document}

\title{Optimal Regularity for Fully Nonlinear Nonlocal Equations with Unbounded Source Terms}

\author{Disson S. dos Prazeres and Makson S. Santos}

\maketitle

\begin{abstract} 
We prove optimal regularity estimates for viscosity solutions to a class of fully nonlinear nonlocal equations with unbounded source terms. More precisely, depending on the integrability of the source term $f \in L^p(B_1)$, we establish that solutions belong to classes ranging from $C^{\sigma-d/p}$ to $C^\sigma$, at critical thresholds. We use approximation techniques and Liouville-type arguments. These results represent a novel contribution, providing the first such estimates in the context of not necessarily concave nonlocal equations.

\medskip

\medskip


\end{abstract}

\maketitle

\tableofcontents



\section{Introduction}

In this article, we examine the regularity theory of viscosity solutions to a fully nonlinear nonlocal equation of the form 

\begin{equation}\label{eq_main}
{\mathcal I}_\sigma(u, x) = C(\sigma)\sup_{\alpha \in {\mathcal A}}\inf_{\beta \in {\mathcal B}}\int_{\mathbb{R}^d}\delta(u, x, y)\dfrac{y^TA_{\alpha,\beta}(x)y}{|y|^{\sigma+d+2}}dy = f(x) \;\;\mbox{ in } \; B_1,   
\end{equation}
where $\sigma \in (0, 2)$, $C(\sigma) > 0$ is a normalizing constant, $\delta(u, x,y):= u(x+y) + u(x-y) -2u(x)$, ${\mathcal A}$ and  ${\mathcal B}$ are indexes sets, and the source term $f$ belongs to a suitable Lesbegue space. Moreover, we suppose the following ellipticity-type condition: for each $\alpha \in {\mathcal A}$ and $\beta \in {\mathcal B}$, $A_{\alpha, \beta}$ is a symmetric $d\times d$ matrix  $(\lambda, \Lambda)$-elliptic, i.e., 
\begin{equation}\label{eq:ellippp}
\lambda I \leq A \leq \Lambda I \quad \mbox{ in }\,\, B_1,    
\end{equation}
for $0 < \lambda \leq \Lambda$. We establish H\"older regularity for solutions (or their gradient), depending on the range of $p$ for which $f \in L^p(B_1)$. 

Regularity results for nonlocal operators have been extensively studied over the years. In \cite{Silvestre2006}, the author presents an analytical proof of Hölder continuity and introduces more flexible assumptions on the operator than previous studies. But it is only in \cite{Caffarelli-Silvestre2009}, where such regularity estimates were given uniformly as the degree of the operator $\sigma \to 2$, seen as a natural extension of second-order operators. In that paper, the authors also prove properties such as comparison principle, a nonlocal Alexandrov-Bakelman-Pucci estimate (ABP for short), Harnack inequality, and regularity in $C^{1, \alpha}$ spaces for viscosity solutions to equations of type \eqref{eq_main}.

Numerous results were uncovered based on the findings in \cite{Caffarelli-Silvestre2009}. In \cite{Caffarelli-Silvestre2011}, the authors extended their previous results using perturbative methods, including the $C^{1, \alpha}$-regularity for a class of non translation-invariant equations. Around the same time, the authors in \cite{Barles-Chasseigne-Imbert-11} also worked with non translation-invariant equations, establishing H\"older regularity of solutions, using different assumptions from \cite{Caffarelli-Silvestre2011}, allowing first-order terms and some degeneracy in the operators. We also mention the work \cite{Kriventsov13}, where the author gives $C^{1, \alpha}$-estimates for more general kernels than previous results.

Besides the $C^{1,\alpha}$ estimates, Evans-Krylov-type results have been explored by the community. In this direction, the authors in \cite{CaffarelliSilvestre2011-2} proved that viscosity solutions for concave (or convex) equations of the form
\begin{equation}\label{eq_conceq}
{\mathcal I}_\sigma(u, x) := \inf_{\alpha \in {\mathcal A}}\int_{\mathbb{R}^d}\delta(u, x, y)K_\alpha(y)dy = 0 \;\;\mbox{ in } \; B_1,   
\end{equation}
are of class $C^{\sigma+\beta}$. They assume that, for each $\alpha \in {\mathcal A}$, the kernel $K_\alpha$ belongs to the class ${\mathcal L}_2$, i.e., they are of class $C^2$ away from the origin, symmetric, satisfies the ellipticity condition
\begin{equation}\label{eq_uniellkernel}
    (2-\sigma)\dfrac{\lambda}{|y|^{d+\sigma}} \leq K_\alpha(y) \leq (2-\sigma)\dfrac{\Lambda}{|y|^{d+\sigma}},
\end{equation}
and 
\begin{equation}\label{eq_l2kernel}
    D^2K_\alpha(y) \leq \dfrac{C}{|y|^{d+2+\sigma}}.
\end{equation}
In \cite{Serra2014Csigma}, the author extended the result above to equations of type \eqref{eq_conceq} with rough kernels, i.e., where the kernels $K_\alpha$ satisfy \eqref{eq_uniellkernel} but not necessarily \eqref{eq_l2kernel}, for every $\alpha \in {\mathcal A}$. Under the additional (and optimal) assumption of $C^\alpha$ exterior data for the solution, the author provides $C^{\sigma+\alpha}$ a priori estimates. Moreover, results involving non translation-invariant kernels are also given, provided the dependency on the variable $x$ is of a $C^\alpha$-fashion. We also mention the work \cite{yu17}, where the author also provides $C^{\sigma+\alpha}$ for viscosity solutions to nonlocal equations of the type
\[
F(D^\sigma u)= 0 \;\;\mbox{ in } \;B_1
\]
provided the quantities $\|u\|_{L^\infty(B_1)}$ and $\|u\|_{L^1_\sigma(\R^d)}$ are sufficiently small.

A natural question is whether one can prove optimal regularity estimates for fully nonlinear nonlocal equations with an unbounded right-hand side. However, such results are relatively scarce in the literature. One of the primary challenges in this context is the lack of compactness properties for viscosity solutions of (1.1) (or even the simpler equation (1.7) below) when $f \in L^p(B_1)$. In addition, stability results are also crucial for implementing the now-classic strategy proposed in \cite{Caffarelli-Silvestre2011}. In fact, only a few works involving an unbounded norm on the right-hand side are available. Moreover, the existing results consider a class of kernels that differ from the broader ${\mathcal L}_0$ class discussed in \cite{Caffarelli-Silvestre2009, Caffarelli-Silvestre2011, Kriventsov13, Serra2014Csigma}, which defines the extremal operator
\[
{\mathcal M}^-_{{\mathcal L}_0}(u,x) = (2-\sigma)\inf_{\lambda \leq a(x,y) \leq \Lambda}\int_{\R^d}\delta(u,x,y)\dfrac{a(x,y)}{|y|^{d+\sigma+2}}.
\]
Instead, they work with kernels of the form
\begin{equation}\label{eq:kernell}
K_{\alpha}(y) = \dfrac{y^TA_{\alpha} y}{|y|^{d+\sigma+2}}dy,    
\end{equation}
where $A_{\alpha}$ is a symmetric $d\times d$ matrix satisfying 
\begin{equation}\label{eq_newellip}
\lambda I \leq \dfrac{1}{d+\sigma}\left(\sigma A_{\alpha} + {\rm Tr}(A_{\alpha})I\right) \leq \Lambda I \quad \mbox{ in }\,\, B_1,    
\end{equation}
for positive constants $0 < \lambda \leq \Lambda$. These kernels define a class of extremal operators as
\[
{\mathcal M}^-(u,x) = \inf_\alpha\int_{\R^d}\delta(u,x,y)K_{\alpha}(y)dy.
\]
It is important to note that, because they allow some degeneracy (conform \eqref{eq_newellip}), although this is a smaller class, it is not necessarily contained within the ${\mathcal M}^-_{{\mathcal L}_0}$.

In this direction, the authors in \cite{GuilenSchwab2012} established a quantitative ABP estimate for viscosity supersolutions to 
\begin{equation}\label{eq_subsolM-}
\left\{
\begin{aligned}
    {\mathcal M}^-(u,x) & \leq f(x) \;\; \mbox{ in } \; B_1 \\
    u(x) & \geq 0 \;\; \mbox{ on } \; \R^d\setminus B_1.
\end{aligned}
\right.
\end{equation}
They show that supersolutions to the equation above satisfy
\[
-\inf_{B_1}u \leq C(d, \lambda)(\|f^+\|_{L^\infty(K_u)})^{(2-\sigma)/2}(\|f^+\|_{L^d(K_u)})^{\sigma/2}, 
\]
where $K_u$ is the coincidence set between $u$ and a type of fractional convex envelope of $u$. In \cite{Kitano20}, the author improved the previous result, by removing the dependency of the $L^\infty$-norm in the estimate above, in case $\sigma$ is sufficiently close to 2.
In \cite{kitano22}, the author removed the restriction on the degree of the operator present in \cite{Kitano20}, by proving that supersolutions to \eqref{eq_subsolM-} with $f \in L^p(B_1)$, for $p \in (d-\varepsilon_0, \infty)$, satisfy
\[
-\inf_{B_1}u \leq C(d, \sigma, \lambda, \Lambda)\|f^+\|_{L^p(B_1)}.
\] 
Moreover, the author also gives $W^{\sigma, p}$-estimates for viscosity solutions to concave equations of the form
\begin{equation}\label{eq_main2}
{\mathcal I}_\sigma(u, x) = C(\sigma)\inf_{\alpha \in {\mathcal A}}\int_{\mathbb{R}^d}\delta(u, x, y)\dfrac{y^TA_{\alpha}(x)y}{|y|^{\sigma+d+2}}dy = f(x) \;\;\mbox{ in } \; B_1,   
\end{equation}
which, in particular, implies that solutions are of class $C^{\sigma-d/p}$. To our knowledge, this is the only result in the literature that gives this type of regularity for fully nonlinear nonlocal equations as in \eqref{eq_main2}, in the presence of an unbounded right-hand side.

The main purpose of this paper is to establish optimal regularity estimates for viscosity solutions to equation \eqref{eq_main}, where the source term $f \in L^p(B_1)$. We emphasize that such estimates have not been previously available for this type of equation. We employ the so-called half-relaxed method, originally introduced in \cite{BarlesPerthame1987}, to prove that a sequence $(u_k)_{k\in \N}$ of merely bounded viscosity solutions to \eqref{eq_main} converges uniformly to a function $u_\infty$ which solves a suitable equation, see Lemma \ref{prop_stability} below. This method relies on the comparison principle of the operator and has been used previously in the context of Hamilton-Jacobi equations and nonlocal equations with Neumann boundary conditions in a half-space. We refer the reader to \cite{barles11, Barles2013}. Once compactness and stability are available, standard approximation arguments can be applied to show the regularity properties of the viscosity solutions.

As discussed earlier, the regularity of solutions depends on the range of $p$ for which $f \in L^p(B_1)$, in the spirit of \cite{Teixeira2014I}. First, we prove that if $p \in (d-\varepsilon_0, \frac{d}{\sigma-1})$, then solutions are of class $C^{\sigma-d/p}$. The constant $\varepsilon_0$ is known as the Escauriaza's exponent.  Notice that this type of regularity was previously known from \cite{kitano22} in the context of a concave equation as in \eqref{eq_main2}. 

The borderline case $f \in L^{d/(\sigma-1)}$ is particularly significant. In the local case, i.e., $\sigma = 2$, this value separates continuity estimates from differentiability properties in the regularity theory. Moreover, this quantity also appears in ABP and Harnack estimates. Meanwhile, in the nonlocal case, this is the first time such a threshold has been explicitly considered, as previous ABP estimates in \cite{GuilenSchwab2012, Kitano20, kitano22} only consider the $L^d$-norm of the right-hand side. We believe this is the correct threshold for future general ABP estimates for the general class ${\mathcal L}_0$. Nevertheless, in this scenario, we show that viscosity solutions to \eqref{eq_main} are Log-Lipschitz, which is better than $C^\alpha$ for every $\alpha \in (0,1)$. 

For $f \in L^p$ for $p \in (\frac{d}{\sigma-1}, +\infty)$, we prove that the solutions belong to the class $C^{1,\alpha}$, where $\alpha$ is defined in \eqref{eq_alpha3}. Finally, for the borderline case where $f \in \rm{BMO}(B_1)$, we show that viscosity solutions to \eqref{eq_main2} are locally of class $C^\sigma$, which is the best regularity we can hope for without assuming further regularity for $f$. These results are detailed in Theorems \ref{main_theo1}, \ref{main_theo2}, \ref{main_theo3} and \ref{main_theo4} below.  



The class of kernels in \eqref{eq_main} follows the form of \eqref{eq:kernell}, as described in \cite{GuilenSchwab2012, Kitano20, kitano22}, but with the additional requirement of uniform ellipticity, as in \eqref{eq:ellippp}. In particular, the class of kernels that we deal with is a subset of ${\mathcal L}_0$ and enjoys all its properties. Furthermore, because uniform ellipticity as in \eqref{eq:ellippp} ensures that the condition \eqref{eq_newellip} is satisfied with the same ellipticity constants, our class of kernels is also included in those described in \cite{GuilenSchwab2012, Kitano20, kitano22}. Developing a similar theory for more general kernels, such as those in ${\mathcal L}_0$, remains an open challenge, primarily due to the need for an appropriate ABP estimate, as discussed earlier.

The remainder of this paper is structured as follows: In section 2 we gather some auxiliary results and present our main results. The proof of the optimal H\"older regularity is the subject of Section 3. In Section 4 we put forward the Log-Lipschitz regularity of solutions. Section 5 is devoted to the proof of H\"older regularity for the gradient of solutions. Finally, in the last section, we investigate the borderline $C^{\sigma}$-regularity.

\section{Preliminaries}

\subsection{Notations and definitions}

This section collects some definitions and notations used throughout the paper. The open ball of radius $r$ and centered at $x_0$ in $\R^d$ is denoted by $B_r(x_0)$. For $\alpha \in (0, 1]$, the notation $u \in C^{\alpha^-}(B_1)$ means that $u \in C^\beta(B_1)$, for every $\beta < \alpha$. We proceed by defining the Log-Lipschitz space.

\begin{definition}
A function $u$ belongs to $C^{\rm Log-Lip}(B_r)$ if there exists a universal constant $C > 0$ such that
\[
\sup_{B_{r/2}(x_0)}|u(x) - u(x_0)| \leq r\ln{r^{-1}}.
\]
\end{definition}

\begin{definition}[BMO space]
We say that $f \in {\rm BMO}(B_1)$ if for all $B_r(x_0) \subset B_1$, we have
\[
\|f\|_{{\rm BMO}(B_1)} := \sup_{0<r\leq1}\dashint_{B_r(x_0)}|f(x) - \langle f\rangle_{x_0,r} |dx, 
\]
where $ \langle f\rangle_{x_0,r} := \displaystyle\dashint_{B_r(x_0)}f(x)dx$.
\end{definition}

Since $u$ is defined in the whole $\R^d$, it can behave very widely as $|x| \to \infty$. Hence, we work within a class where we have a certain decay of the solutions as they approach infinity, see also \cite{Caffarelli-Silvestre2009}. 

\begin{definition}[Growth at infinity]
We say that a function $u:\mathbb{R}^d \rightarrow \mathbb{R}$ belongs to $L^1_\sigma(\R^d)$, if
\[
\|u\|_{L^1_\sigma(\R^d)}:= \int_{\R^d} |u(x)|\dfrac{1}{1 +|x|^{d+\sigma}}dx < +\infty.
\]
\end{definition}

In the next, we define viscosity solutions:

\begin{definition}[Viscosity solution]
We say that an upper (lower) semicontinuous function $u \in L_\sigma^1(\R^d)$ is a viscosity subsolution (supersolution) to \eqref{eq_main}, if for any $x_0 \in B_1$ and $\varphi \in C^2(B_r(x))$ such that $u-\varphi$ has a local maximum (minimum) at $x_0$, then the function
\[
v(x):=
\left\{
\begin{aligned}
    \varphi(x) &\;\; \mbox{ in } \; B_r(x_0) \\
    u(x) &\;\; \mbox{ in } \; \R^d\setminus B_r(x_0),
\end{aligned}
\right.
\]
satisfies
\[
-{\mathcal I}_\sigma(v,x) \leq (\geq) f(x) \;\; \mbox{ in } B_1
\]
We say a function $u \in C(\overline{B}_1)\cap L_\sigma^1(\R^d)$ is a viscosity solution to \eqref{eq_main} if it is simultaneously a viscosity subsolution and supersolution.
\end{definition}

\begin{remark}[Scaling properties]\label{rem_scaling}
Throughout the manuscript, we assume certain smallness conditions on the norms of $u$ and the source term $f$. We want to stress that such conditions are not restrictive. In fact, if $u$ is a viscosity solution to \eqref{eq_main}, then for $\varepsilon >0$ the function 
\[
v(x) := \dfrac{u(x)}{\|u\|_{L^\infty(B_1)} + \varepsilon^{-1}\|f\|_{L^p(B_1)}},
\]
satisfy $\|v\|_{L^\infty(B_1)} \leq 1$ and solves
\[
{\mathcal I}_\sigma(v, x) = \tilde{f}(x) \;\; \mbox{ in } \; B_1, 
\]
where 
\[
\tilde{f}(x) := \dfrac{f(x)}{\|u\|_{L^\infty(B_1)} + \varepsilon^{-1}\|f\|_{L^p(B_1)}},
\] 
is such that $\|\tilde{f}\|_{L^\infty(B_1)} \leq \varepsilon$.
\end{remark}

\subsection{Main results}

As mentioned earlier, the regularity of solutions depends on the range of $p$ under consideration. The critical cases occur when $f \in L^{\frac{d}{\sigma-1}}(B_1)$ for $C^\alpha$-regularity and $f \in \mathrm{BMO}$ for $C^{\sigma}$-regularity. Due to the nonlocal nature of the problem, achieving these critical cases requires $f$ to possess higher regularity compared to the local case. For instance, in the local case, the first critical threshold is at $p = d$, leading to $C^{\text{Log-Lip}}_{\text{loc}}$-regularity. Similarly, when $f \in \mathrm{BMO}$ in the local case, it yields $C^{1,\text{Log-Lip}}_{\text{loc}}$-regularity.

We now present the main results of this article, beginning with a result in Hölder spaces for the case where $p$ is below $\frac{d}{\sigma-1}$.

\begin{theorem}\label{main_theo1}
Let $u \in C(B_1)\cap L^1_\sigma(\R^d)$ be a viscosity solution to \eqref{eq_main}, with $f \in L^p(B_1)$, $p \in \left(d-\varepsilon_0, \frac{d}{\sigma-1}\right)$. Then, $u \in C^\alpha_{loc}(B_1)$ for any
\begin{equation}\label{eq_alpha1}
\alpha \in \left(0, \frac{\sigma p -d}{p}\right].    
\end{equation}
Moreover, there exists a positive constant $C=C(p,d, \sigma, \lambda, \Lambda)$, such that
\[
\|u\|_{C^{\alpha}(B_{1/2})} \leq C(\|u\|_{L^\infty(B_1)} + \|u\|_{L^1_\sigma(\R^d)} + \|f\|_{L^p(B_1)}). 
\]
\end{theorem}

We observe that for $\sigma = 2$, we have $p \in (d - \varepsilon_0, d)$ and $\alpha \in (0, \frac{2p-d}{p})$, recovering the regularity result for the local case  reported in \cite{Teixeira2014I}. Next, we consider the borderline case $p = \frac{d}{\sigma - 1}$. In this scenario, we show that solutions are Log-Lipschitz continuous, achieving the same level of regularity as in \cite[Theorem 2]{Teixeira2014I}, but with the requirement of higher regularity for the source term $f$.

\begin{theorem}\label{main_theo2}
Let $u \in C(B_1)\cap L^1_\sigma(\R^d)$ be a viscosity solution to \eqref{eq_main}, with $f \in L^p(B_1)$, $p = \frac{d}{\sigma-1}$. Then, $u \in C^{Log-Lip}_{loc}(B_1)$, and there exists a positive constant $C=C(d, \sigma, \lambda, \Lambda)$, such that
\[
\|u\|_{C^{Log-Lip}(B_{1/2})} \leq C\left(\|u\|_{L^\infty(B_1)} + \|u\|_{L^1_\sigma(\R^d)} + \|f\|_{L^{\frac{d}{\sigma-1}}(B_1)}\right). 
\]
\end{theorem}

In what follows, we present our third main theorem. As before, we recover the local regularity in the limit as $\sigma \to 2$, demonstrated in \cite[Theorem 3]{Teixeira2014I}. Recall that $\alpha_0$ comes from the $C^{1, \alpha_0}$-regularity of ${\mathcal I}_\sigma$-harmonic functions. 

\begin{theorem}\label{main_theo3}
Let $u \in C(B_1)\cap L^1_\sigma(\R^d)$ be a viscosity solution to \eqref{eq_main}, with $f \in L^p(B_1)$, $p \in \left(\frac{d}{\sigma-1}, +\infty\right)$. Then, $u \in C^{1, \alpha}_{loc}(B_1)$ for any
\begin{equation}\label{eq_alpha3}
\alpha \in \left(0, \min\left(\sigma - 1 - \frac{d}{p}, \alpha_0^- \right)\right].    
\end{equation}
Moreover, there exists a positive constant $C=C(p, d, \sigma, \lambda, \Lambda)$, such that
\[
\|u\|_{C^{1, \alpha}(B_{1/2})} \leq C(\|u\|_{L^\infty(B_1)} + \|u\|_{L^1_\sigma(\R^d)} + \|f\|_{L^p(B_1)}). 
\]
\end{theorem}

Observe that as $p \to \infty$, the corresponding Hölder exponent approaches $\sigma - 1$, indicating that we achieve $C^\sigma$-regularity in the case $p = \infty$, which is indeed the case as confirmed in \cite{vivas}. However, we also show that this result holds under the weaker assumption that $f \in \mathrm{BMO}(B_1)$, which is a proper subset of $L^\infty(B_1)$. This is the content of our last main result.

\begin{theorem}\label{main_theo4}
Let $u \in C(B_1)\cap L^1_\sigma(\R^d)$ be a viscosity solution to \eqref{eq_main2}, with $f \in {\rm BMO}(B_1)$. Then, $u \in C^{\sigma}_{loc}(B_1)$ and there exists a positive constant $C=C(d, \sigma, \lambda, \Lambda)$, such that
\[
\|u\|_{C^{\sigma}(B_{1/2})} \leq C(\|u\|_{L^\infty(B_1)} + \|u\|_{L^1_\sigma(\R^d)} + \|f\|_{{\rm BMO}(B_1)}). 
\]
\end{theorem}

The proof of Theorem \ref{main_theo4} differs from the strategies used for the previous theorems. This is primarily because scaling of the form $x \mapsto \rho^{-\sigma}u(\rho x)$, for $\rho \ll 1$, leads to a growth rate of $|x|^{-\sigma}$ at infinity, which increases too rapidly to be integrable with respect to the tails of our kernel. Consequently, we employ techniques similar to those in \cite{Serra2014Csigma, Serra2015}, where Liouville-type results are used to establish interior regularity of solutions through blow-up arguments.



\subsection{Auxiliary results}

In this subsection, we prove some results used throughout the paper. Since we did not find any references stating exactly what we needed, we start with a comparison principle for viscosity solutions of \eqref{eq_main} (and also \eqref{eq_main2}), in the case where $f \equiv 0$. See also \cite[Theorem 5.2]{Caffarelli-Silvestre2009}. 

\begin{proposition}[Comparison principle]\label{comparison principle}
Let $u,v \in L^1_\sigma(\R^d)$, $u$ upper semicontinuous and $v$ lower semicontinuous, be respectively a viscosity subsolution and supersolution to the equation
\begin{equation}\label{eq comparison}
{\mathcal I}_\sigma(w, x) = 0 \;\; \mbox{ in } \; B_1,
\end{equation}
such that $u \leq v$ in $\mathbb{R}^N\setminus B_1$. Then, $u \leq v$ in $B_1$.
\end{proposition} 
\begin{proof}
Suppose by contradiction that 
\[
\theta=\sup_{B_1}(u-v)>0.
\]
For $\varepsilon>0$ we define the auxiliary function
\[
\Phi(x,y)=u(x)-v(y)-\frac{|x-y|^2}{2\varepsilon^2},
\]
and consider $(x_\varepsilon,y_\varepsilon)\in \overline{B}_1\times\overline{B}_1$ such that
\[
\Phi(x_\varepsilon,y_\varepsilon)=\sup_{x,y\in B_1}\Phi(x,y).
\]
Notice that $\Phi(x_\varepsilon,y_\varepsilon)\geq \sup_{x\in B_1}\Phi(x,x)=\theta$, which yields to
\begin{equation}\label{inequality 1 cp}
\frac{|x_\varepsilon-y_\varepsilon|}{2\varepsilon^2}\leq u(x_\varepsilon)-v(y_\varepsilon)-\theta.
\end{equation}
By compactness, we have that $(x_\varepsilon,y_\varepsilon) \rightarrow (\bar{x},\bar{y})\in \overline{B}_1\times\overline{B}_1$ as $\varepsilon \to 0$, and by using \eqref{inequality 1 cp} we obtain $\bar{x}=\bar{y}$. Therefore
\[
0 \leq \lim_{\varepsilon\rightarrow 0}\frac{|x_\varepsilon-y_\varepsilon|^2}{2\varepsilon^2} = u(\bar{x}) - v(\bar{x}) - \theta \leq 0, 
\]
which implies
\[
u(\bar{x}) - v(\bar{x}) = \theta > 0. 
\]
Moreover, since $u\leq v$ on $ \R^d\setminus B_1$, we have $\bar{x}\in B_1$. 
We set $\varphi_1(x):= v(y_\varepsilon) + |x-y_\varepsilon|/2\varepsilon^2$ and $\varphi_2(y):= u(x_\varepsilon) + |y-x_\varepsilon|/2\varepsilon^2$, and observe that $u-\varphi_1$ has a local maximum at $x_\varepsilon$, while $v-\varphi_2$ has a local minimum at $y_\varepsilon$. Hence by using that $u$ is a subsolution and $v$ is a supersolution of \eqref{eq comparison}, we have the viscosity inequalities
\begin{equation*}
\sup_{\alpha \in {\mathcal A}}\inf_{\beta \in {\mathcal B}}\left(\int_{B_r}\delta(\varphi_1, x_\varepsilon, y)\dfrac{y^TA_{\alpha,\beta}(x_\varepsilon)y}{|y|^{\sigma+d+2}}dy  + \int_{\R^d\setminus B_r}\delta(u, x_\varepsilon, y)\dfrac{y^TA_{\alpha,\beta}(x_\varepsilon)y}{|y|^{\sigma+d+2}}dy \right) \geq 0,
\end{equation*}
and 
\begin{equation*}
\sup_{\alpha \in {\mathcal A}}\inf_{\beta \in {\mathcal B}}\left(\int_{B_r}\delta(\varphi_2, y_\varepsilon, y)\dfrac{y^TA_{\alpha,\beta}(y_\varepsilon)y}{|y|^{\sigma+d+2}}dy  + \int_{\R^d\setminus B_r}\delta(v, y_\varepsilon, y)\dfrac{y^TA_{\alpha,\beta}(y_\varepsilon)y}{|y|^{\sigma+d+2}}dy \right) \leq 0,
\end{equation*}
for $r$ small enough. Hence, from the definition of $\sup$ and $\inf$, there exist $\alpha \in {\mathcal A}$ and $\beta \in {\mathcal B}$ such that
\begin{equation*}
\varepsilon^{-2} o_r(1) + \int_{\R^d\setminus B_r}\delta(u, x_\varepsilon, y)\dfrac{y^TA_{\alpha,\beta}(x_\varepsilon)y}{|y|^{\sigma+d+2}}dy \geq -\gamma/2,
\end{equation*}
and
\begin{equation*}
\varepsilon^{-2} o_r(1) + \int_{\R^d\setminus B_r}\delta(v, y_\varepsilon, y)\dfrac{y^TA_{\alpha,\beta}(y_\varepsilon)y}{|y|^{\sigma+d+2}}dy \leq \gamma/2,    
\end{equation*}
for $\gamma> 0$ sufficiently small. Subtracting the inequalities above yields to
\[
\varepsilon^{-2} o_r(1) + \int_{\R^d\setminus B_r}\delta(u, x_\varepsilon, y)\dfrac{y^TA_{\alpha,\beta}(x_\varepsilon)y}{|y|^{\sigma+d+2}}dy - \int_{\R^d\setminus B_r}\delta(v, y_\varepsilon, y)\dfrac{y^TA_{\alpha,\beta}(y_\varepsilon)y}{|y|^{\sigma+d+2}}dy \geq -\gamma.
\]
Notice that, When $\varepsilon\rightarrow 0$, by the contradiction hypotheses we have
$$
 \int_{B_1\setminus B_r}\delta(u, \bar{x}, y)\dfrac{y^TA_{\alpha,\beta}(\bar{x})y}{|y|^{\sigma+d+2}}dy-\int_{B_1\setminus B_r}\delta(v, \bar{x}, y)\dfrac{y^TA_{\alpha,\beta}(\bar{x})y}{|y|^{\sigma+d+2}}dy\leq 0.
$$
Moreover,
$$
 \int_{\R^d\setminus B_1}\delta(u, \bar{x}, y)\dfrac{y^TA_{\alpha,\beta}(\bar{x})y}{|y|^{\sigma+d+2}}dy-\int_{\R^d\setminus B_1}\delta(v, \bar{x}, y)\dfrac{y^TA_{\alpha,\beta}(\bar{x})y}{|y|^{\sigma+d+2}}dy\leq - \lambda \theta \int_{\mathbb{R}^d\setminus B_1 } |z|^{-(d + \sigma)}dz.
$$
Therefore
\begin{equation*}
-\lambda \theta \int_{\mathbb{R}^d\setminus B_1 } |z|^{-(d + \sigma)}dz \geq -\gamma,
\end{equation*}
which is a contradiction for $\gamma$ small enough. This finishes the proof.
\end{proof}

We now focus on one of the main contributions of this paper: the stability of solutions to \eqref{eq_main} when $f \in L^p(B_1)$. To the best of our knowledge, this is the first time such a result has been established in the fully nonlinear nonlocal context. For comparison, see \cite[Corollary 4.7]{Caffarelli-Silvestre2009}.

\begin{proposition}[Compactness and stability]\label{prop_stability}
Let $u_k$ be a normalized viscosity solution to
\begin{equation*}
{\mathcal I}_\sigma(u_k, x) = f_k \;\;\mbox{ in } \; B_1.
\end{equation*} 
Suppose that there exists a positive constant $M$ such that 
\begin{equation}\label{eq_condcres1}
    |u_k(x)| \leq M(1+|x|)^{1+\alpha} \;\;\mbox{ for all } \; x \in \R^d.
\end{equation}
Suppose further that
\[
\|f_k\|_{L^p(B_1)} \to 0,
\]
for $p\in (d-\varepsilon_0,+\infty)$. Then there exists $u_\infty \in C(B_1)\cap L_\sigma^1(\R^d)$ such that 
\[
\|u_k-u_\infty\|_{L^\infty(B_{4/5})} \to 0.
\]
Moreover, $u_\infty$ solves
\begin{equation}\label{eq_01}
{\mathcal I}_\sigma(u_\infty, x) = 0 \;\;\mbox{ in } \; B_1.    
\end{equation}
\end{proposition}

\begin{proof}
First, observe that given $R > 0$, we obtain from \eqref{eq_condcres1} 
\[
|u_k(x)| \leq C(R) \;\; \mbox{ for any } \; x \in B_R.
\]
Hence, Given any compact set $\Omega \subset \R^d$, we have that the a.e. limits
\[
\bar{u}(x) := \limsup_{k\to\infty, y_k\to x} u_k(y_k), \;\; x \mbox{ in } \; \Omega,  
\]
and
\[
\underline{u}(x) := \liminf_{k\to\infty, y_k\to x} u_k(y_k), \;\; x \mbox{ in } \; \Omega, 
\]
are well-defined. Since the a.e. convergence holds for every compact set of $\R^d$, we also have the a.e. convergence in the whole $\R^d$. Using this fact and once again \eqref{eq_condcres1}, the Dominated Convergence Theorem ensures that
\begin{equation}\label{eq_07}
\|\bar{u}-u_k\|_{L^1_\sigma(\R^d)} \to 0, 
\end{equation}
and 
\begin{equation}\label{eq_08}
\|\underline{u}-u_k\|_{L^1_\sigma(\R^d)} \to 0,   
\end{equation}
through the respective subsequences.
We are going to show that $\bar{u}$ is a viscosity subsolution to \eqref{eq_01}, and 
$\underline{u}$ is a viscosity supersolution to \eqref{eq_01}. We will prove the subsolution case since the supersolution case is analogous. Let $x_0 \in B_1$ and $\varphi \in C^2(B_r(x_0))$ be such that $\bar{u} - \varphi$ has a maximum at $x_0$. Without loss of generality, we can assume that $\varphi$ is defined by
\begin{equation*}
\varphi = \left\{
\begin{aligned}
P, &\quad\mbox{in}\quad B_r(x_0),\\
\bar{u},&\quad \mbox{on}\quad B_r(x_0)^c,
\end{aligned}
\right.
\end{equation*}
for some paraboloid $P$. We need to show that 
\begin{equation}\label{eq_sub1}
{\mathcal I}_\sigma(\varphi, x) \geq 0 \;\;\mbox{ in }\; B_1.    
\end{equation}
Suppose by contradiction that 
\begin{equation}\label{eq_02}
{\mathcal I}_\sigma(\varphi, x) < -\eta,    
\end{equation}
for some $\eta > 0$. Now, let $\psi_k$ be a viscosity solution to
\begin{equation}\label{eq_04}
\left\{
\begin{aligned}
{\mathcal M}^+_{\lambda^*, 1}(\psi_k, x)  & = & -|f_k(x)| & \;\; \mbox{ in } \; B_r(x_0) \\
     \psi_k & = & 0 & \;\; \mbox{ on } \; \partial B_r(x_0),    
\end{aligned}
\right.    
\end{equation}
for some $\lambda^* < 1$ to be chosen later. Here, the maximal operator ${\mathcal M}^+$ is defined with respect to the class ${\mathcal L}_0$, as in \cite{Caffarelli-Silvestre2009} (and defined below). We have 
\begin{align*}
{\mathcal I}_\sigma(\varphi + \psi_k,x) - {\mathcal I}_\sigma(\varphi, x) & \leq {\mathcal M}^+_{\lambda, \Lambda}(\psi_k, x) \\
& = \Lambda\int_{\mathbb{R}^n}\dfrac{\delta^+(\psi_k, x, y)}{|y|^{n+\sigma}}dy - \lambda\int_{\mathbb{R}^n}\dfrac{\delta^-(\psi_k, x, y)}{|y|^{n+\sigma}}dy \\
& = \Lambda\int_{\mathbb{R}^n}\dfrac{\delta^+(\psi_k, x, y)}{|y|^{n+\sigma}}dy - \dfrac{\lambda}{\lambda^*}\left(\int_{\mathbb{R}^n}\dfrac{\delta^+(\psi_k, x, y)}{|y|^{n+\sigma}}dy + |f_k(x)|\right) \\
& = \left(\Lambda - \dfrac{\lambda}{\lambda^*} \right)\int_{\mathbb{R}^n}\dfrac{\delta^+(\psi_k, x, y)}{|y|^{n+\sigma}}dy - \dfrac{\lambda}{\lambda^*}|f_k(x)|,
\end{align*}
where we have used \eqref{eq_04} to conclude
\[
\int_{\mathbb{R}^n}\dfrac{\delta^+(\psi_k, x, y)}{|y|^{n+\sigma}}dy - \lambda^*\int_{\mathbb{R}^n}\dfrac{\delta^-(\psi_k, x, y)}{|y|^{n+\sigma}}dy = -|f_k(x)|.
\]
Now, by choosing $\lambda^*$ such that $\Lambda - \dfrac{\lambda}{\lambda^*} \leq 0$ and $-\dfrac{\lambda}{\lambda^*} \leq -1$, we obtain
\begin{equation}\label{eq_03}
{\mathcal I}_\sigma(\varphi + \psi_k,x) \leq {\mathcal I}_\sigma(\varphi, x) + f_k(x).    
\end{equation}
Let $P_k$ be defined by
\begin{equation}
		\varphi_k=\left\{
		\begin{aligned}
			P, &\quad\mbox{in}\quad B_r(x_0),\\
			u_k,&\quad \mbox{on}\quad B_r(x_0)^c,
		\end{aligned}
		\right.
	\end{equation}
and 
\begin{equation}
		\Tilde{P}=\left\{
		\begin{aligned}
			c|x-x_0|^2, &\quad\mbox{in}\quad B_r(x_0),\\
			0,&\quad \mbox{on}\quad B_r(x_0)^c.
		\end{aligned}
		\right.
	\end{equation}
By using the ABP estimates in \cite[Theorem 3.1]{kitano22} we have $\|\psi_k\|_\infty\rightarrow 0$, as $k\rightarrow \infty$. Then, there exists a $x_k\in B_r(x_0) $ such that $\varphi_k+\psi_k+\Tilde{P}$ touch $u_k$ by above in $B_r(x_0)$. Therefore, we have the viscosity inequality
\begin{equation}\label{contradction fx_k}
{\mathcal I}_\sigma(\varphi_k+\psi_k+\Tilde{P})\geq f(x_k).
\end{equation}
By ellipticity and using \eqref{eq_02} and \eqref{eq_03} we obtain that
\begin{equation}\label{eq_06}
{\mathcal I}_\sigma(\varphi_k+\psi_k+\Tilde{P}) \leq {\mathcal I}_\sigma(\varphi_k+\psi_k)+{\mathcal M}^+_\sigma(\Tilde{P})
 \leq {\mathcal I}_\sigma(\varphi_k+\psi_k)-{\mathcal I}_\sigma(\varphi+\psi_k)-\eta+ f(x_k)+{\mathcal M}^+_\sigma(\Tilde{P}).  
\end{equation}
We observe that
\begin{align*}
\Big|\int_{\R^N}\delta(\varphi_k+\psi_k,x,y)dy  -\int_{\R^N}\delta(\varphi+\psi_k,x,y)dy\Big| & \leq \int_{\R^N\setminus B_r(x_0)}|\delta(u_k,x,y)-\delta(\bar{u},x,y)|dy \\
& \leq C(r)\|u_k-\bar{u}\|_{L^1_\sigma(\R^d)},    
\end{align*}
and by using \eqref{eq_07}, we can conclude that for $k$ sufficiently large
\begin{equation}\label{eq_05}
|{\mathcal I}_\sigma(\varphi_k+\psi_k)-{\mathcal I}_\sigma(\varphi+\psi_k)| \leq \dfrac{\eta}{4}.
\end{equation}
Hence, from \eqref{eq_06} and \eqref{eq_05}, we obtain
\[
{\mathcal I}_\sigma(\varphi_k+\psi_k+\Tilde{P})\leq \dfrac{\eta}{4}  -\eta +f(x_k)+{\mathcal M}^+_\sigma(\Tilde{P}).
\]
Now, choose $c(\Lambda,N,\sigma)$ sufficiently small so that
\[
{\mathcal M}^+_\sigma(\Tilde{P})\leq \eta/4.
\]
Finally, for $k$ sufficiently large we get
\[
{\mathcal I}_\sigma(\varphi_k+\psi_k+\Tilde{P},x_k)\leq -\eta/2 +f(x_k),
\]
which is a contradiction with \eqref{contradction fx_k}. This finishes the proof of \eqref{eq_sub1}, i.e., $\bar{u}$ solves in the viscosity sense
\[
{\mathcal I}_\sigma(\bar{u}, x) \geq 0 \;\;\mbox{ in }\; B_1, 
\]
We similarly show that
\begin{equation*}
    {\mathcal I}_\sigma(\underline{u}, x) \leq 0 \;\;\mbox{ in }\; B_1.  
\end{equation*}
Now, from the definition of $\bar{u}$, $\underline{u}$ and the viscosity inequalities above we can infer from Proposition \ref{comparison principle}, that in fact 
\[
\bar{u} = \underline{u} = u_\infty,
\]
and hence up to a subsequence, $u_k \to u_\infty$ locally uniformly in $B_1$ (see for instance \cite[Lemma 6.2]{Barles2013}). Moreover, from the viscosity inequalities satisfied by $\bar{u}$ and $\underline{u}$, we conclude that $u_\infty$ solves
\[
{\mathcal I}_\sigma(u_\infty, x) = 0 \;\;\mbox{ in }\; B_1,
\]
in the viscosity sense.
\end{proof}

Using the stability result above, we can prove the following Approximation Lemma, which relates the solutions to our problem with ${\mathcal I}_\sigma$-harmonic functions. 

\begin{lemma}[Approximation Lemma]\label{lem_approximation}
Let $u \in C(B_1)\cap L^1_\sigma(\R^d)$ be a normalized viscosity solution to \eqref{eq_main}, with $p \in (d-\varepsilon_0, +\infty)$. Suppose that 
\[
|u(x)|\leq M(1+|x|)^{1+\alpha} \;\;\mbox{ for all }\;\; x\in \R^d.
\]
Given $\delta >0$ there exist $\varepsilon >0$, such that if
$$
\|f\|_{L^p(B_1)}\leq \varepsilon,
$$
we can find a function $h\in C^{1,\alpha_0}(B_{4/5})$ satisfying 
\[
\sup_{B_{3/4}}|u-h|\leq \delta.
\]
\end{lemma}

\begin{proof}
Suppose not, then there exist $\delta_0 > 0$ and sequences $(u_k)_{k\in\mathbb N}$, $(f_k)_{k\in\mathbb N}$ such that 
\begin{equation}\label{eq_10}
{\mathcal I}_\sigma(u_k, x) = f_k \;\;\mbox{ in } \; B_1,\\
\end{equation}  
\begin{equation}\label{eq_11}
    \|f_k\|_{L^p(B_1)} \leq \dfrac{1}{k},
\end{equation}
and 
\begin{equation}\label{eq_12}
    |u_k(x)|\leq M(1+|x|)^{1+\alpha} \;\;\mbox{ for all }\;\; x\in \R^d,
\end{equation}
but,
\begin{equation}\label{eq_12.5}
    |u_k-h|>\delta_0,
\end{equation}
for all $h\in C^{1,\alpha_0}(B_{4/5})$. From the contradiction hypotheses \eqref{eq_10}, \eqref{eq_11}, \eqref{eq_12} and Proposition \ref{prop_stability}, we can guarantee the existence of a function $u_\infty \in C(B_1)\cap L^1_\sigma(\R^d)$ such that $u_k \to u_\infty$ locally uniformly in $B_1$ satisfying
\begin{equation}\label{eq_13}
{\mathcal I}_\sigma(u_\infty, x) = 0 \;\;\mbox{ in } \; B_1.\\
\end{equation}  
Now, the regularity available for \eqref{eq_13}, see \cite{Caffarelli-Silvestre2009}, implies that $u_\infty \in C^{1,\alpha_0}(B_{4/5})$. By taking $h\equiv u_\infty$, we reach a contradiction with \eqref{eq_12.5} for $k$ sufficiently large.
\end{proof}

\section{H\"{o}lder regularity}

In this section, we detail the proof of Theorem \ref{main_theo1}, namely, the optimal $C^\alpha_{loc}$-regularity, for
\[
\alpha \in \left(0, \frac{\sigma p -d}{p}\right], 
\]
where $p \in \left(d-\varepsilon_0, \frac{d}{\sigma-1}\right)$. We start by applying Lemma \ref{lem_approximation} and showing the existence of a constant close to $u$ in sufficiently small balls. 

\begin{proposition}\label{lem:step1Calpha}
Let $u \in C(B_1)\cap L^1_\sigma(\R^d)$ be a normalized a viscosity solution to \eqref{eq_main}, with $p \in \left(d-\varepsilon_0, \frac{d}{\sigma-1}\right)$.
Assume that 
\[
|u(x)|\leq M(1+|x|)^{1+\alpha} \;\;\mbox{ for all }\;\; x\in \R^d.
\]
If
\[
\|f\|_{L^p(B_1)}\leq \varepsilon,
\]
then there exist constants $0 < \rho \ll 1/2$ and $A$ satisfying, $|A| \leq C$ and  
\[
\sup_{B_\rho}|u-A| \leq \rho^\alpha,
\]
where $C> 0$ is a universal constant.
\end{proposition}

\begin{proof}
Fix $\delta > 0$ (to be chosen later) and let $h$ be the function from Lemma \ref{lem_approximation}. Since $h \in C^{1, \alpha_0}(B_{4/5})$, for $\rho$ sufficiently small, we have
\[
|h - h(0)| \leq C\rho. 
\]
Now, from Lemma \ref{lem_approximation} and the Triangular inequality we obtain
\begin{align*}
\sup_{B_\rho}|u - h(0)| & \leq  \sup_{B_\rho}|u - h| + \sup_{B_\rho}|h - h(0)|  \\
& \leq \delta + C\rho.
\end{align*}
Now, we make the universal choices 
\begin{equation}\label{eq_choice}
\rho = \min\left[\left(\frac{1}{2C}\right)^{\frac{1}{1-\alpha}}, \left(\frac{1}{(1+C)100}\right)^{\frac{1}{\alpha_0-\alpha}}\right]  \;\;\mbox{ and } \;\; \delta = \frac{\rho^\alpha}{2},    
\end{equation}
and by setting $A = h(0)$, we conclude that
\[
\sup_{B_\rho}|u - A| \leq \rho^\alpha. 
\]
Notice that the choice of $\delta$ determines the value of $\varepsilon$ via Lemma \ref{lem_approximation}.
\end{proof}

In what follows, we iterate the previous proposition to find a sequence of constants that approaches $u$ at the origin. 

\begin{proposition}\label{lem:step2Calpha}
Let $u \in C(B_1)\cap L^1_\sigma(\R^d)$ be a normalized a viscosity solution to \eqref{eq_main}, with $p \in \left(d-\varepsilon_0, \frac{d}{\sigma-1}\right)$. Assume that
\[
|u(x)|\leq M(1+|x|)^{1+\alpha} \;\;\mbox{ for all }\;\; x\in \R^d.
\]
If
\[
\|f\|_{L^p(B_1)}\leq \varepsilon,
\]
then we can find a sequence $(A_k)_{k \in \mathbb{N}}$ satisfying
\begin{equation}\label{eq_21}
\sup_{B_{\rho^k}}|u(x)-A_k| \leq \rho^{k\alpha},    
\end{equation}
with
\begin{equation}\label{eq_20}
|A_{k+1} - A_k| \leq C\rho^{k\alpha}.     
\end{equation}    
\end{proposition}

\begin{proof}
We argue by an induction argument. By setting $A_0=A_1 = 0$, the case $k=0$ follows immediately. Suppose we have verified the statement for $k = 1, \ldots, n$, and let us prove the case $k=n+1$. We introduce the auxiliary function $v_k : \R^d \longrightarrow \R$ 
\[
v_k(x) := \frac{u(\rho^kx) - A_k}{\rho^{k\alpha}}.
\]
Notice that by \eqref{eq_21} we have $|v_k(x)| \leq 1$ in $B_1$. In addition, $v_k$ solves
\[
\mathcal{I}_\sigma(v,x) = \tilde{f}(x)\;\; \mbox{ in } \; B_1,
\]
where $\tilde{f}(x) = \rho^{k(\sigma-\alpha)}f(\rho^kx)$. Moreover, our choice of $\alpha$ in \eqref{eq_alpha1} assures $\|\tilde{f}\|_{L^p(B_1)} \leq \varepsilon$. Next, we are going to show that $v_k$ satisfies
\begin{equation}\label{eq_22}
 |v(x)| \leq M(1+|x|^{1+\alpha_0}) \;\; \mbox{ for all } \;\; x \in \R^d,   
\end{equation}
for some universal constant $M$. In fact, we resort again to an induction argument. For $k=0$, we have $v_0=u$ and \eqref{eq_22} is verified. Now, assume that the case $k=1, \ldots, n$ is already verified. We shall prove the case $k=n+1$. Observe that
\[
v_{n+1} = \dfrac{v_{n}(\rho x) - \tilde{A}_n}{\rho^\alpha},
\]
where $\tilde{A}_n$ comes from Lemma \eqref{lem_approximation} applied to $v_n$. Now, for $2|x|\rho > 1$ we estimate
\begin{align*}
|v_{n+1}(x)| &\leq \rho^{-\alpha} \left(|v_n(\rho x)| +|\tilde{A}_n|\right)\\
&\leq \rho^{-(1+\alpha)}\left[(1+ \rho^{1+\alpha_0}|x|^{1+\alpha_0}) + C(1 + \rho|x|)\right]\\
&\leq  \rho^{(\alpha_0 - \alpha)}(5 +9C)|x|^{1 + \alpha_0}\\
&\leq  |x|^{1+ \alpha_0},    
\end{align*}
where in the last inequality we used \eqref{eq_choice}. On the other hand, if $2|x|\rho \leq 1$, we obtain
\begin{align*}
|v_{n+1}(x)| & \leq \rho^{-\alpha}(|v_n(\rho x) - \tilde{h}(\rho x)| + |\tilde{h}(\rho x) -\tilde{A}_n| \\
&\leq \rho^{-\alpha}\left(\dfrac{\rho^{\alpha}}{2}+ C\rho|x|\right) \\
&\leq \frac{1}{2} + \frac{C}{2\rho^{\alpha}}\\
&\leq M(1+|x|^{1+\alpha_0}),    
\end{align*}
where $M:= 1/2 + C/(2\rho^{\alpha})$, and hence \eqref{eq_22} is proved. Finally, we now can apply Proposition \ref{lem:step1Calpha} to $v_k$ and we obtain
\[
\sup_{B_\rho}|v_k - \tilde{A}_k| \leq \rho^\alpha,
\]
and rescaling back to $u$ we conclude
\[
\sup_{B_{\rho^{k+1}}}|u - A_{k+1}| \leq \rho^{(k+1)\alpha},
\]
where $A_{k+1} = A_k + \rho^{k\alpha}\tilde{A}_k$, which satisfies \eqref{eq_20}. This finishes the proof.
\end{proof}

We are now ready to prove Theorem \ref{main_theo1}.

\begin{proof}[Proof of Theorem \ref{main_theo1}]
Notice that from \eqref{eq_20} we have that $(A_k)_{k\in\N}$ is a Cauchy sequence, and hence there exists $A_\infty$ such that $A_k \to A_\infty$, as $k\to \infty$. Moreover, we also have from \eqref{eq_20}
\[
|A_k-A_\infty| \leq C\rho^{k\alpha}.
\]
Now, fix $0< r \ll 1$ and let $k \in \N$ be such that $\rho^{k+1} \leq r \leq \rho^k$. We estimate,
\begin{align*}
\sup_{B_r}|u(x)-A_\infty| & \leq \sup_{B_{\rho^k}}|u(x)-A_k| + \sup_{B_{\rho^k}}|A_k - A_\infty| \\
& \leq \rho^{k\alpha} + C\rho^{k\alpha} \\
& \leq \frac{(C+1)}{\rho^{\alpha}}\rho^{(k+1)\alpha} \\
& \leq Cr^\alpha.
\end{align*}
By taking the limit as $k \to \infty$ in \eqref{eq_21}, we obtain $A_\infty = u(0)$. This finishes the proof.
\end{proof}

\section{Log-Lipschitz continuity}

This section addresses the first critical case $p = \frac{d}{\sigma - 1}$, which yields the desired Log-Lipschitz regularity. In particular, solutions are of class $C_{loc}^\alpha$ for every $\alpha \in (0, 1)$. As before, we begin by demonstrating the existence of a linear approximation of $u$ within sufficiently small balls.

\begin{proposition}
Let $u \in C(B_1)\cap L^1_{\sigma}(\R^d)$ be a normalized viscosity solution to \eqref{eq_main}, with $p =\frac{d}{\sigma-1}$. Suppose further that
\[
|u(x)| \leq M(1+|x|)^{1+\alpha_0} \;\; \mbox{ for all } \;\; x\in \R^d.
\]
If
\[
\|f\|_{L^p(B_1)} \leq \varepsilon,
\]
then, there exist a constant $0 < \rho \ll 1/2$ and an affine function $\ell$ of the form
\[
\ell(x) = A + B\cdot x,
\]
satisfying $|A|, |B| \leq C$ and 
\[
\sup_{B_\rho}|u(x)-\ell(x)| \leq \rho. 
\]    
\end{proposition}

\begin{proof}
The proof is similar to Proposition \ref{lem:step1Calpha}. We fix $\delta >0$ to be determined later. For $\rho \ll 1/2$, we have that
\[
|h(x)-h(0)-Dh(0)\cdot x| \leq C\rho^{1+\alpha_0},
\]
where $h \in C^{1,\alpha_0}(B_{3/4})$ comes from Lemma \ref{lem_approximation}. By setting $\ell(x) = h(0) + Dh(0)\cdot x$, we obtain from the Triangular inequality that
\begin{align*}
\sup_{B_\rho}|u(x)-\ell(x)| & \leq \sup_{B_\rho}|u(x)-h(x)| + \sup_{B_\rho}|h(x)-\ell(x)| \\ 
& \leq \delta + C\rho^{1+\alpha_0}.
\end{align*}
As before, we make universal choices
\begin{equation}\label{eq_choices2}
\rho = \left(\frac{1}{2C}\right)^{\frac{1}{\alpha_0}} \;\; \mbox{ and } \;\; \delta = \frac{\rho}{2},   
\end{equation}
which determines the value of $\varepsilon$ through Lemma \ref{lem_approximation}. Therefore
\[
\sup_{B_\rho}|u(x)-\ell(x)| \leq \rho.
\]
\end{proof}

\begin{proposition}\label{prop:step2C1alpha}
Let $u \in C(B_1)\cap L^1_{\sigma}(\R^d)$ be a normalized viscosity solution to \eqref{eq_main}, with $p=\frac{d}{\sigma-1}$. Suppose further that
\[
|u(x)| \leq M(1+|x|)^{1+\alpha_0} \;\; \mbox{ for all } \;\; x\in \R^d.
\]
If
\[
\|f\|_{L^p(B_1)} \leq \varepsilon,
\]
then, there exists a sequence of affine function $(\ell_k)_{k \in \N}$ of the form
\[
\ell_k(x) = A_k + B_k\cdot x,
\]
satisfying
\begin{equation}\label{condition AB loglipschtz}
\frac{|A_{k+1}-A_k|}{\rho^{k}} + |B_{k+1}-B_k| \leq C,     
\end{equation}
and
\[
\sup_{B_\rho}|u(x)-\ell(x)| \leq \rho^{k}. 
\]    
\end{proposition}
\begin{proof}
As before, we resort to an induction argument. By considering $\ell_0=\ell_1=0$, the case $k=0$ follows trivially. Now, suppose that the cases $k=1, \ldots, n$ have been verified, and let us prove the case $k=n+1$. We define the auxiliary function $v_k: \R^d \longrightarrow \R$ by
\[
v_k(x) := \frac{u(\rho^kx)-\ell(\rho^kx)}{\rho^{k}}. 
\]
We have that $v_k(x) \leq 1$ in $B_1$, and solves
\[
{\mathcal I_\sigma}(x,v) = \tilde{f}(x) \;\; \mbox{ in } \;\; B_1, 
\]
where $\tilde{f}(x) = \rho^{k(\sigma-1)}f(\rho^kx)$. Notice that as $\sigma-1>0$, we have $\|\tilde{f}\|_{L^p(B_1)} \leq \varepsilon$. Arguing similarly as in Proposition \ref{lem:step2Calpha}, we can also show that 
\[
|v_k(x)| \leq 1 + |x|^{1+\alpha_0}.
\]
Hence, we can apply Proposition \ref{prop:step1C1alpha} to $v_k$ to conclude that there exists $\tilde{\ell}_k=\tilde{A}+\tilde{B}x$ such that
\[
\sup_{B_\rho}|v_k(x)-\tilde{\ell}_k(x)| \leq \rho.
\]
Rescaling back to $u$, we obtain
\[
\sup_{B_\rho}|u(x)-\ell_{k+1}(x)| \leq \rho^k,
\]
where $\ell_{k+1} = \ell_k(x) + \rho^{k}\tilde{\ell}_k(\rho^{-1}x)$. Observe that 

$$
|A_{k+1}-A_k|=|\tilde{A}\rho^k|\leq C\rho^k,
$$
and
$$
|B_{k+1}-B_k|=|\tilde{B}|\leq C
$$
which prove the condition \eqref{condition AB loglipschtz}.
\end{proof}

We now present the proof of Theorem \ref{main_theo2}

\begin{proof}[Proof of Theorem \ref{main_theo2}]
Notice that, by \eqref{condition AB loglipschtz}, $(A_k)_{k\in\N}$ is a Cauchy sequence, and hence, there exists $A_\infty$ such that $A_k \to A_\infty$, as $k\to \infty$. Moreover, we have
\[
|A_k-A_\infty| \leq C\rho^{k\alpha}.
\]
Now, fix $0< r \ll 1$ and let $k \in \N$ be such that $\rho^{k+1} \leq r \leq \rho^k$. We have,
\begin{align*}
\sup_{B_r}|u(x)-A_\infty| & \leq \sup_{B_{\rho^k}}|u(x)-A_k-B_kx| + \sup_{B_{\rho^k}}|B_kx| \\
& \leq \rho^{k} + Ck\rho^k \\
& \leq \frac{1}{\rho}\left(\rho^{k+1}+Ck\rho^{k+1}\right) \\
& \leq Cr+\dfrac{\ln{r}}{\ln{\rho}}Cr\\
&\leq -Cr\ln{r}.
\end{align*}
Finally, by taking the limit as $k \to \infty$ in \eqref{eq_21}, we obtain $A_\infty = u(0)$. This finishes the proof.
\end{proof}

\section{H\"{o}der continuity of the gradient}

In this section, we give the proof of Theorem \ref{main_theo3}, in which we prove $C_{loc}^{1, \alpha}$-regularity for
\begin{equation*}
\alpha \in \left(0, \min\left(\sigma-1-\frac{d}{p}\,, \alpha_0^-\right)\right].    
\end{equation*}
The proof follows the general lines of the proof of Theorem \ref{main_theo1}, but now at the gradient level.

\begin{proposition}\label{prop:step1C1alpha}
Let $u \in C(B_1)\cap L^1_{\sigma}(\R^d)$ be a normalized viscosity solution to \eqref{eq_main}, with $p \in \left(\frac{d}{\sigma-1}, +\infty\right)$. Suppose further that
\[
|u(x)| \leq M(1+|x|)^{1+\alpha_0} \;\; \mbox{ for all } \;\; x\in \R^d.
\]
If
\[
\|f\|_{L^p(B_1)} \leq \varepsilon,
\]
then, there exist a constant $0 < \rho \ll 1/2$ and an affine function $\ell$ of the form
\[
\ell(x) = A + B\cdot x,
\]
satisfying $|A|, |B| \leq C$ and 
\[
\sup_{B_\rho}|u(x)-\ell(x)| \leq \rho^{1+\alpha}. 
\]    
\end{proposition}

\begin{proof}
We fix $\delta >0$ to be determined later. For $\rho \ll 1/2$, we have that
\[
|h(x)-h(0)-Dh(0)\cdot x| \leq C\rho^{1+\alpha_0},
\]
where $h \in C^{1,\alpha_0}(B_{3/4})$ comes from Lemma \ref{lem_approximation}. By setting $\ell(x) = h(0) + Dh(0)\cdot x$, we obtain from the Triangular inequality that
\begin{align*}
\sup_{B_\rho}|u(x)-\ell(x)| & \leq \sup_{B_\rho}|u(x)-h(x)| + \sup_{B_\rho}|h(x)-\ell(x)| \\ 
& \leq \delta + C\rho^{1+\alpha_0}.
\end{align*}
As before, we make universal choices
\begin{equation}\label{eq_choices2}
\rho = \min\left[\left(\frac{1}{2C}\right)^{\frac{1}{\alpha_0-\alpha}}, \left(\frac{1}{(1+C)100} \right)^{\frac{1}{\alpha_0-\alpha}} \right] \;\; \mbox{ and } \;\; \delta = \frac{\rho^{1+\alpha}}{2},   
\end{equation}
which determines the value of $\varepsilon$ through Lemma \ref{lem_approximation}. Therefore
\[
\sup_{B_\rho}|u(x)-\ell(x)| \leq \rho^{1+\alpha}.
\]
\end{proof}

\begin{proposition}\label{prop:step2C1alpha}
Let $u \in C(B_1)\cap L^1_{\sigma}(\R^d)$ be a normalized viscosity solution to \eqref{eq_main}, with $p \in \left(\frac{d}{\sigma-1}, +\infty\right)$. Suppose further that
\[
|u(x)| \leq M(1+|x|)^{1+\alpha_0} \;\; \mbox{ for all } \;\; x\in \R^d.
\]
If
\[
\|f\|_{L^p(B_1)} \leq \varepsilon,
\]
then, there exists a sequence of affine function $(\ell_k)_{k \in \N}$ of the form
\[
\ell_k(x) = A_k + B_k\cdot x,
\]
satisfying
\begin{equation}\label{eq_condAB}
|A_{k+1}-A_k| + \rho^{k}|B_{k+1}-B_k| \leq C\rho^{k(1+\alpha)},     
\end{equation}
and
\[
\sup_{B_{\rho^k}}|u(x)-\ell_k(x)| \leq \rho^{k(1+\alpha)}. 
\]    
\end{proposition}

\begin{proof}
As before, we resort to an induction argument. By considering $\ell_0=\ell_1=0$, the case $k=0$ follows trivially. Now, suppose that the cases $k=1, \ldots, n$ have been verified, and let us prove the case $k=n+1$. We define the auxiliary function $v_k: \R^d \longrightarrow \R$ by
\[
v_k(x) := \frac{u(\rho^kx)-\ell(\rho^kx)}{\rho^{k(1+\alpha)}}. 
\]
We have that $v_k(x) \leq 1$ in $B_1$, and solves
\[
{\mathcal I_\sigma}(x,v) = \tilde{f}(x) \;\; \mbox{ in } \;\; B_1, 
\]
where $\tilde{f}(x) = \rho^{k(\sigma-\alpha-1)}f(\rho^kx)$. Notice that for our choice of $\alpha$ in \eqref{eq_alpha3}, we have $\|\tilde{f}\|_{L^p(B_1)} \leq \varepsilon$. Arguing similarly as in Proposition \ref{lem:step2Calpha}, we can also show that 
\[
|v_k(x)| \leq 1 + |x|^{1+\alpha_0}.
\]
Hence, we can apply Proposition \ref{prop:step1C1alpha} to $v_k$ to conclude that there exists $\tilde{\ell}_k$ such that
\[
\sup_{B_\rho}|v_k(x)-\tilde{\ell}_k(x)| \leq \rho^{1+\alpha}.
\]
Rescaling back to $u$, we obtain
\[
\sup_{B_\rho}|u(x)-\ell_{k+1}(x)| \leq \rho^{1+\alpha},
\]
where $\ell_{k+1} = \ell_k(x) + \rho^{k(1+\alpha)}\tilde{\ell}_k(\rho^{-1}x)$. From the definition of $\ell_{k+1}$, it is immediate that condition \eqref{eq_condAB} is also satisfied. This finishes the proof.
\end{proof}

\begin{proof}[Proof of Theorem \ref{main_theo3}]
The proof follows the same lines as in Theorem 1.1 and Theorem 1.2. Notice that from \eqref{eq_condAB} we have that $(A_k)_{k\in\N}$ and $(B_k)_{k\in\N}$ are a Cauchy sequence, and hence, we can find $A_\infty$ and $B_\infty$ satisfying $A_k \to A_\infty$ and $B_k \to B_\infty$, as $k\to \infty$. Moreover, we have
\[
|A_k-A_\infty| \leq C\rho^{k(1+\alpha)} \;\;\;\mbox{ and }\;\;\; |B_k-B_\infty| \leq C\rho^{k\alpha}.
\]
Now, fix $0< r \ll 1$ and let $k \in \N$ be such that $\rho^{k+1} \leq r \leq \rho^k$. We have,
\begin{align*}
\sup_{B_r}|u(x)-A_\infty-B_\infty\cdot x| & \leq \sup_{B_{\rho^k}}|u(x)-A_k-B_k\cdot x| + \sup_{B_{\rho^k}}|A_k - A_\infty| + \rho^k\sup_{B_{\rho^k}}|B_k - B_\infty| \\
& \leq \rho^{k(1+\alpha)} + C\rho^{k(1+\alpha)} + C\rho^{k(1+\alpha)}\\
& \leq \frac{C}{\rho^{1+\alpha}}\rho^{(k+1)(1+\alpha)}\\
& \leq Cr^{1+\alpha}.
\end{align*}
Finally, by taking the limit as $k \to \infty$ in \eqref{eq_21}, we obtain $A_\infty = u(0)$. We can also show that $B_\infty = Du(0)$, see for instance \cite{Bronzi-Pimentel-Rampasso-Teixeira2020}. This finishes the proof.
\end{proof}

\section{The borderline case}

This section deals with the $C^\sigma$-regularity for viscosity solutions of \eqref{eq_main2}. As previously discussed, since $1+\alpha = \sigma$, we can no longer follow the strategy employed above. In this case, we follow the ideas put forward in \cite{Serra2014Csigma} (see also\cite{yu17}). We begin with a technical lemma, that can be found in \cite[Claim 3.2]{Serra2014Csigma}. 

\begin{lemma}\label{caracterizacao de c1alpha}
Let $0<\overline{\alpha}<\alpha<1$ and $u\in C^{1,\overline{\alpha}}(B_1)$. If $u(0)=|Du(0)|=0$ and
$$
\sup_{0<r<1/2}r^{\overline{\alpha}-\alpha}[u]_{C^{1,\overline{\alpha}}(B_r)}\leq A,
$$
for a constant $A>0$, then
\[
[u]_{C^{1+\alpha}(B_{1/2})} \leq 2A.
\]
\end{lemma}

\begin{proposition}\label{prop:Csigma}
Let $u \in C^{1+\alpha}(B_1)\cap L^1_{\sigma}(\R^d)$, with $\alpha< \sigma-1$, be a normalized viscosity solution to \eqref{eq_main2} with $f\in {\rm BMO}(B_1)$. Suppose further that
\[
\|u\|_{C^{1+\alpha}(B_1)} \leq M \quad \mbox{ and }\quad \|f\|_{{\rm BMO}(B_1)} \leq \varepsilon,
\]
then $u\in C^\sigma(B_{1/2})$ and
\[
\|u\|_{C^{\sigma}(B_{1/2})} \leq C.
\]
\end{proposition}

\begin{proof}
We argue by contradiction. Suppose that the result is false, then we can find sequences $(u_k)_{k\in\mathbb{N}}, (f_k)_{k\in\mathbb{N}}$, such that for every $k\in \mathbb{N}$, 
\begin{align}
\begin{split}\label{eq:contprop6,2}
[u_k]_{C^{1+\alpha}(B_1)} \leq M \\
\|f_k\|_{{\rm BMO}(B_1)} < \varepsilon_k,
\end{split}    
\end{align}
but
\[
\|u_k\|_{C^{\sigma}(B_{1/2})} > k,
\]
where $\varepsilon_k \to 0$, as $k \to \infty$. For $k \in \N$, we define the quantity 
$$
   \theta_k(r')=\sup_{r'<r<1/2}\sup_{z\in B_{1/2}}r^{1+\alpha - \sigma}[u]_{C^{1+\alpha}(B_r(z))}
$$
and observe that
$$
\lim_{r'\rightarrow 0}\theta_k(r')=\sup_{r'>0}\theta_k(r').
$$
Moreover, if $r_1 \leq r_2$, then $\theta_k(r_2) \leq \theta_k(r_1)$. By the Lemma \ref{caracterizacao de c1alpha} we have that
$$
\sup_{r'>0}\theta_k(r')\geq k/2,
$$
therefore there exists a $r_k>1/k$ and $z_k\in B_{1/2}$ such that
\begin{equation}\label{contradição BMO}
r_k^{1+\alpha - \sigma}[u_k]_{C^{1+\alpha}(B_{r_k}(z_k))}>\theta_k(1/k)>\theta_k(r_k)\rightarrow \infty.
\end{equation}
Since from \eqref{eq:contprop6,2} $[u_k]_{C^{1+\alpha}(B_1)}<M$, we have $r_k\rightarrow 0$ as $k \to \infty$. Now, for $R \in \left[1, \frac{1}{2r_k}\right]$, we define the blow-up $v_k:B_R\rightarrow \mathbb{R}$
\[
v_k(x)=\frac{1}{\theta_k(r_k)}\frac{1}{r_k^{\sigma}}u_k(r_kx+z_k).
\]
Notice that
\begin{align}
    \begin{split}\label{estimativa para vk}
     [v_k]_{ C^{1+\alpha}(B_{R})} & =\frac {r_k^{1+\alpha-\sigma}}{\theta_k(r_k)}[u_k]_{C^{1+\alpha}(B_{r_kR}(z_k))} \\
     &  = \frac{(Rr_k)^{1+\alpha-\sigma}[u_k]_{B_{r_kR}(z_k)}}{\theta_k(r_k)}R^{\sigma-1-\alpha}\\
     & \leq R^{\sigma-1-\alpha}.   
    \end{split}
\end{align}
Consider the auxiliary function $w_k: \R^d \to \R$ defined as
\[
w_k(x) := (v_k-l_k)(x),
\]
where  
\[
l_k=v_k(0)+Dv_k(0)x.
\]
Since that $[l_k]_{ C^{1+\alpha}(B_{R})}=0$, it follows from (\ref{estimativa para vk}) that
\begin{equation}\label{estimativa de wk}
[w_k]_{ C^{1+\alpha}(B_{R})}\leq R^{\sigma-1-\alpha}.
\end{equation}
In particular, for $R=1$ we have
$$
[Dw_k]_{ C^{\alpha}(B_{1})}\leq 1,
$$
therefore
$$
|Dw_k(x)|=|Dw_k(x)-Dw_k(0)|\leq |x|^{\alpha}, \quad \mbox{ for }\quad x\in B_1,
$$
which implies
\begin{equation}\label{eq:estB1}
|w_k(x)|=|w_k(x)-w_k(0)|\leq \|Dw_k\|_{L^\infty(B_{|x|})}|x|\leq |x|^{1+\alpha},    
\end{equation}
for all $x \in B_1$. Let $\eta$ be a smooth function such that $\eta=1$ in $B_{1/2}$ and $\eta=0$ outside $B_1$. For $e\in \mathbb{S}^n$ we have
$$
\int_{B_1} \eta\cdot D_e w_k dx=\int_{B_1}D_e \eta\cdot w_k dx\leq C(n).
$$
Hence, there exists $z\in B_1$ such that $|Dw_k(z)|\leq C(n)$ and by (\ref{estimativa de wk}) we have
$$
|D_ew_k(x)-D_ew_k(z)|\leq R^{\sigma-1-\alpha}|x-z|^{\alpha}.
$$
Therefore, for $x\in B_R$ and $1\leq R \leq \frac{1}{2r_k}$,
\begin{equation}\label{eq_helpconbeta}
|D_ew_k(x)|\leq C(n)+ R^{\sigma-1-\alpha}|x-z|^{\alpha}\leq CR^{\sigma-1}.    
\end{equation}
Our goal now is to show that
\begin{equation}\label{eq_goalbeta}
[w_k]_{C^\beta(B_R)}\leq R^{\sigma-\beta},    
\end{equation}
for all $\beta \in [0,1+\alpha]$ and $R \in \left[1, \frac{1}{2r_k}\right]$.

{\it Case $\beta =0$:} Notice that for $1\leq |x|\leq R$, we have from \ref{eq_helpconbeta}
\begin{align}
\begin{split}\label{eq:beta00}
|w_k(x)| & = |w_k(x)-w_k(0)| \\
& \leq \|D_ew_k(x)\|_{L^\infty(B_{|x|})}|x| \\
& \leq CR^{\sigma}.
\end{split}
\end{align}
The estimate for $x \in B_1$ follows from \eqref{eq:estB1}.

{\it Case $\beta \in (0, 1)$:} In this case, we estimate for $x,\bar{x}\in B_R$
\begin{align*}
    |w_k(x)-w_k(\bar{x})| &\leq \|Dw_k\|_{L^{\infty}(B_R)}|x-\bar{x}| \\
    & \leq R^{\sigma-1}|x-\bar{x}|^{1-\beta}|x-\bar{x}|^\beta \\
    & \leq CR^{\sigma-\beta}|x-\bar{x}|^\beta,
\end{align*}
which gives
\begin{equation}\label{eq:beta01}
[w_k]_{ C^{\beta}(B_R)}\leq CR^{\sigma-\beta}.    
\end{equation}

{\it Case $\beta \in [1, 1+\alpha]$:} Finally, we have
\begin{align*}
   |Dw_k(x)-Dw_k(\bar{x})| & \leq [Dw_k]_{C^{\alpha}(R)}|x-\bar{x}|^{\alpha}\\
   & \leq CR^{\sigma-1-\alpha}R^{\alpha-\beta+1}|x-\bar{x}|^{\beta-1},
\end{align*}
which implies
$$
[Dw_k]_{ C^{\beta-1}(B_R)}\leq CR^{\sigma-\beta},
$$
or equivalently
\begin{equation}\label{eq:b1alpha}
[w_k]_{ C^{\beta}(B_R)}\leq CR^{\sigma-\beta}.
\end{equation}
Therefore, \eqref{eq_goalbeta} follows from \eqref{eq:beta00}, \eqref{eq:beta01} and \eqref{eq:b1alpha}. Thus, there exists a $w\in C^{1+\alpha}(\mathbb{R}^d)$ such that $ w_k\rightarrow w$ locally uniformly in the $C^{1+\alpha}$-norm. Now, from \eqref{contradição BMO}, we have that
$$
[w_k]_{C^{1+\alpha}(B_1)}\geq \frac{1}{2}
$$
and therefore
\begin{equation}\label{eq_contnorm}
[w]_{C^{1+\alpha}(B_1)}\geq \frac{1}{2}.    
\end{equation}
Moreover, from \eqref{eq_goalbeta},
\begin{equation}\label{eq_goalbetaw}
[w]_{C^\beta(\R^d)}\leq R^{\sigma-\beta},    
\end{equation}
for all $\beta \in [0,1+\alpha]$ and $R > 1$. Notice that $w_k$ solves 
\[
{\mathcal I}_\sigma(w_k, x) = \tilde{f}_k(x),
\]
where $\tilde{f}_k(x) = \dfrac{1}{\theta_k(r_k)}f_k(r_kx+z_k)$,
for $p\geq d$, we have
\begin{eqnarray*}
\|\tilde{f}_k(x)\|_{L^p(B_1}    &= &   \left ( \int_{B_{1}} |\tilde{f}_k(x)|^p dx \right )^{1/p}  \\
 &= & |B_1|^{1/p} \left ( \intav{B_{r_k}(z_k)} |f_k(x)|^p dx \right )^{1/p}  \\
 &\le & |B_1|^{1/p}\|f_k\|_{\rm BMO(B_1)} \le |B_1|^{1/p} \varepsilon_k.
\end{eqnarray*}
Hence,
\begin{align*}
    {\mathcal M}^+_{{\mathcal L}_0}(w_k(\cdot + h) - w_k) & \geq {\mathcal I}_\sigma(w_k, x+h) - {\mathcal I}_\sigma(w_k, x) \\
    & = \dfrac{1}{\theta_k(r_k)}\left[ f_k(r_k(x+h) + z_k) - f_k(r_kx + z_k) \right] .
\end{align*}
Moreover,
$$
|w_k(x + h) - w_k|\leq [w_k]_{ C^{\alpha}(B_R)}|h|^\alpha\leq C|x|^{\sigma-\alpha}
$$
 and as $[f_k]_{\rm BMO(B_1)}\rightarrow 0$ (and therefore $\|f_k\|_{L^p(B_1)} \to 0$), we have from Proposition \ref{prop_stability} that
\[
{\mathcal M}^+_{{\mathcal L}_0}(w(\cdot + h) - w) \geq 0 \;\;\ \mbox{ in } \;\; \R^d.
\]
Similarly, we can prove that
\[
    {\mathcal M}^-_{{\mathcal L}_0}(w(\cdot + h) - w) \leq 0\;\;\ \mbox{ in } \;\; \R^d.
\]
Therefore, 
\begin{equation}\label{eq:m+-0}
    {\mathcal M}^-_{{\mathcal L}_0}(w(\cdot + h) - w) \leq 0 \leq {\mathcal M}^+_{{\mathcal L}_0}(w(\cdot + h) - w) \;\;\ \mbox{ in } \;\; \R^d.
\end{equation}
Finally, for
$$
\tilde{w}_k=\dashint w_k(\cdot + h)d\mu(h) - w_k,
$$
we have that
$$
|\tilde{w}_k|\leq \dashint |w_k(\cdot + h) - w_k|d\mu(h)\leq C|x|^{\sigma-1-\alpha}.
$$
The concavity of ${\mathcal I}_\sigma$ yields 
\begin{align*}
{\mathcal M}^+_{{\mathcal L}_0}\left( \dashint\tilde w_k(\cdot + h)d\mu(h) - w_k , x\right) & \geq {\mathcal I}_\sigma\left(\dashint\tilde w_k(\cdot + h)d\mu(h), x\right) - {\mathcal I}_\sigma(w_k, x) \\
& \geq \dashint{\mathcal I}_{\sigma}(w_k, x+h)d\mu(h) - \tilde{f}_k(x) \\
& = \dashint \tilde{f}_k(x+h) - \tilde{f}_k(x) d\mu(h).
\end{align*}
Hence, by passing the limit as $k \to \infty$ we get
\begin{equation}\label{eq:m+0av}
    {\mathcal M}^+_{{\mathcal L}_0}\left( \dashint w(\cdot + h)d\mu(h) - w , x \right) \geq 0.
\end{equation}
Therefore, from \eqref{eq_goalbetaw}, \eqref{eq:m+-0} and \eqref{eq:m+0av}, we can apply Theorem \cite[Theorem 2.1]{Serra2014Csigma}, to conclude that $w$ is a polynomial of degree 1, which is a contradiction with \eqref{eq_contnorm}.

\end{proof}

\begin{proof}[Proof of Theorem \ref{main_theo4}] Recall that by Remark \ref{rem_scaling} we can assume $\|f\|_{BMO(B_1)} \leq \varepsilon$, where  $\varepsilon>0$ is the constant from the previous lemma. Since $\|f\|_{L^p(B_1)} \le C\|f\|_{\rm BMO(B_1)}$, we can use Theorem \ref{main_theo3} to conclude 
\[
\|u\|_{C^{1, \alpha}(B_{1/2})} \leq C(\|u\|_{L^\infty(B_1)} + \|u\|_{L^1_\sigma(\R^d)} + \|f\|_{{\rm BMO}(B_1)}).
\]
Therefore from Proposition \ref{prop:Csigma} we have
\[
\|u\|_{C^{\sigma}(B_{1/2})} \leq C(\|u\|_{L^\infty(B_1)} + \|u\|_{L^1_\sigma(\R^d)} + \|f\|_{BMO(B_1)}),
\]
proving the result.
\end{proof}

\bigskip
{\bf Acknowledgement:} D. dos Prazeres was partially supported by CNPq and CAPES/Fapitec. M. Santos was partially supported by the Portuguese government through FCT-Funda\c c\~ao para a Ci\^encia e a Tecnologia, I.P., under the projects UID/MAT/04459/2020, and PTDC/MAT-PUR/1788/2020. This study was financed in part by the Coordena\c c\~ao de Aperfei\c coamento de Pessoal de N\'ivel Superior - Brazil (CAPES) - Finance Code 001.

\bibliographystyle{plain}
\bibliography{prazeres_santos}

\bigskip

\noindent\textsc{Makson S. Santos}\\
Departamento de Matem\'atica do Instituto Superior T\'ecnico\\
Universidade de Lisboa\\
1049-001 Lisboa, Portugal\\
\noindent\texttt{makson.santos@tecnico.ulisboa.pt}

\vspace{.15in}

\noindent\textsc{Disson dos Prazeres}\\
Department of Mathematics\\
Universidade Federal de Sergipe - UFS, \\
49100-000, Jardim Rosa Elze, São Cristóvão - SE, Brazil\\
\noindent\texttt{disson@mat.ufs.br}

\end{document}